\documentclass{amsart}
\usepackage{amssymb,latexsym,amsmath,amsthm}
\newtheorem{theorem}{Theorem}
\newtheorem*{thm*}{Theorem}
\newtheorem{prop}{Proposition}
\newtheorem*{prop*}{Proposition}
\newtheorem{lemma}[prop]{Lemma}
\newtheorem*{lemma*}{Lemma}

\newtheorem*{cor*}{Corollary}
\newtheorem*{conjecture*}{Conjecture}

\theoremstyle{definition}

\newcommand{\bstack}[2]{\substack{#1\\#2}}

\newcommand{\maps}{\rightarrow}

\newcommand{\al}{\alpha}
\newcommand{\be}{\beta}

\newcommand{\ga}{\gamma}
\newcommand{\del}{\delta}
\newcommand{\ep}{\epsilon}

\newcommand{\om}{\omega}
\newcommand{\Om}{\Omega}

\newcommand{\lam}{\lambda}

\newcommand{\Ga}{\Gamma}

\newcommand{\Rcal}{\mathcal{R}}

\newcommand{\Tcal}{\mathcal{T}}

\newcommand{\C}{\mathbb{C}}

\newcommand{\HH}{\mathbb{H}}

\newcommand{\R}{\mathbb{R}}

\newcommand{\Z}{\mathbb{Z}}

\newcommand{\x}{{\bf x}}
\newcommand{\y}{{\bf y}}

\newcommand{\beq}{\begin{equation}}
\newcommand{\eeq}{\end{equation}}

\begin{document}

\title[Discrete fractional integral operators]{A note on discrete fractional integral operators on the Heisenberg group}

\author{Lillian B. Pierce}
\address{School of Mathematics, Institute for Advanced Study, Princeton New Jersey 08540}
\email{lbpierce@math.ias.edu}
\thanks{The author is supported by the Simonyi Fund at the Institute for Advanced Study and the National Science Foundation, including DMS-0902658 and DMS-0635607.}
\date{22 April 2010.}

\maketitle
\begin{abstract}
We consider the discrete analogue of a fractional integral operator on the Heisenberg group, for which we are able to prove nearly sharp results by means of a simple argument of a combinatorial nature. 
\end{abstract}

\section{Introduction}
In the setting of Euclidean space, singular Radon transforms take the form
\beq\label{R_dfn}
\Rcal f(x) = p. v. \int_{\R^{k_0}} f(\ga_t(x)) K(t)dt,
\eeq
where $K$ is a Calder\'{o}n-Zygmund kernel and $\ga_t(\cdot)$ is a family of diffeomorphisms of $\R^k$ that depend smoothly on $t \in \R^{k_0}$. It is also required that $\ga_0$ be the identity mapping and the family of varieties $\{ \ga_t(x): t \in \R^{k_0} \}_{x \in \R^k}$ have certain curvature properties generalizing the property of finite type; these curvature properties may be stated either in terms of iterates of the mappings $\ga_t$ or commutation properties of associated vector fields. Such Radon transforms are now well understood; the general theory may be found in \cite{CNSW}, which also points to earlier literature in the field.
Less well understood are fractional Radon transforms:  for example operators of the form (\ref{R_dfn}) but with the Calder\'{o}n-Zygmund kernel $K(t)$ replaced by $|t|^{-k_0\lam}$, where $0 < \lam < 1$.  It is not yet understood in general how the $(L^p, L^q)$ bounds of such operators depend on the geometry of the underlying family of submanifolds. 

One specific fractional Radon transform for which sharp results are known is the fractional integral operator on the Heisenberg group $\HH^k$, defined by
\beq\label{T_HH_dfn} \mathcal{T}^\lam f (z,t) = f * (|w|^{-2k\lam} \cdot \del_{\tau=0})(z,t) = \int_{\C^k} f(z-w,t- \langle z,w \rangle) |w|^{-2k \lam} dw.
\eeq
Here $0< \lam < 1$, $(z,t)$ and $(\om,\tau)$ are elements in $\HH^k$ (identified topologically with $\C ^k \times \R$), $\del$ denotes the Dirac delta function, and $\langle z,w \rangle = 2 \Im(z \cdot \bar{w})$ is the symplectic bilinear form associated with the Heisenberg group. 
Work of Ricci and Stein \cite{RSIII} and Christ \cite{Chr88} provides the following sharp result for $\mathcal{T}^\lam$: 
\begin{theorem}[Ricci and Stein; Christ]\label{Christ_thm}
Let $\Gamma_k$ denote the closed triangle in $[0,1]^2$ with vertices  $(0,0), (1,1),$ and $(\frac{2k+1}{2k+2}, \frac{1}{2k+2})$. The operator $\Tcal^\lam$ extends to a bounded operator from $L^p(\HH^k)$ to $L^q(\HH^k)$  if and only if $(1/p,1/q) \in \Ga_k$, $ \lam<1$, and $p,q$ satisfy the homogeneity condition
\[ \frac{1}{q} = \frac{1}{p} - \frac{2k(1-\lam)}{2k+2}. \]
\end{theorem}
(Here we have left the factor of 2 in the homogeneity condition in order to highlight the role played by the homogeneous dimension $2k+2$ of $\HH^k$.)

Recently, attention has turned to discrete analogues of operators of the form (\ref{R_dfn}), (\ref{T_HH_dfn}), and their relatives. In this note we consider the discrete analogue of $\Tcal^\lam$, acting on (compactly supported) functions $f : \Z^{2k+1} \maps \C$ by
\beq\label{T_dis}
T^\lam f(n,t) = \sum_{\bstack{m \in \Z^{2k}}{m \neq 0}} \frac{f(n-m,t-\om( n,m))}{|m|^{2k\lam}},
\eeq
where $0< \lam <1$, $n = (n_1,n_2), m=(m_1,m_2) \in \Z^k \times \Z^k, t \in \Z$, and  $\om$ is the symplectic bilinear form defined by $\om(n,m) = 2(n_2 \cdot m_1 - n_1 \cdot m_2)$, in analogy to the form $\langle \cdot, \cdot \rangle$ on $\HH^k$. 
We prove
a nearly sharp (within $\ep$) result for $T^\lam$ in dimension $k=1$. 
\begin{theorem}\label{frxnl_int_thm_dim1}
For $k=1$ and $0< \lam < 1$, $T^\lam$ extends to a bounded operator from $\ell^p(\Z^{3})$ to $\ell^q(\Z^{3})$ with
\[ ||T^\lam f||_{\ell^q(\Z^3)} \leq A ||f||_{\ell^p(\Z^3)}\]
 if $p,q$ satisfy 
\begin{itemize}
\item[(i)]  $1/q < 1/p - \frac{1}{2} (1-\lam),$
\item[(ii)] $1/q< \lam, 1/p> 1-\lam$.
\end{itemize}
\end{theorem}
Theorem \ref{frxnl_int_thm_dim1} would be sharp if it included $p,q$ satisfying equality in condition (i). Conditions (i) (including equality) and (ii) are in fact necessary for $T^\lam$ to be a bounded operator from $\ell^p$ to $\ell^q$; we provide the relevant examples at the end of this paper.

To situate this result in the context of current knowledge of discrete analogues, we note certain key results by way of background. In discrete analogues of Radon transforms, the curvature conditions for $\ga_t(\cdot)$ in (\ref{R_dfn}) are replaced by polynomial structure. Thus the discrete analogue of (\ref{R_dfn}) takes the form 
\[ Rf(n) = \sum_{\bstack{m \in \Z^{k_0}}{m \neq 0}} f(P(n,m))K(m),\]
where $P$ is a polynomial mapping from $\Z^{k} \times \Z^{k_0}$ to $\Z^k$ and $K$ is a Calder\'{o}n-Zygmund kernel on $\R^{k_0}$. The only progress made to date for such discrete singular 
Radon transforms is in the translation invariant and quasi-translation invariant cases. 

In the translation invariant case, in which $P(n,m) = n-Q(m)$ with $Q$ a polynomial mapping from $\Z^{k_0}$ to $\Z^k$, Ionescu and Wainger \cite{IW} (building on earlier work of \cite{AO}, \cite{SW0}) have proved the deep result that $R$ is bounded on $\ell^p(\Z^k)$ for all $1< p< \infty$, with $\ell^p$ norm dependent only on the degree of the polynomial $Q$ and independent of its coefficients. 
Similarly, there has been significant recent progress on translation invariant discrete fractional Radon transforms, in which $K(m)$ is replaced by $|m|^{-k_0\lam}$ for any $0<\lam <1$ (see  \cite{SW2}, \cite{SW3}, \cite{Obe}, \cite{Pie_J_lam}, \cite{Pie_Waring}).

In the quasi-translation invariant case, a discrete singular Radon transform takes the form
\beq\label{quasi_R}
Rf(n, n') = \sum_{\bstack{m \in \Z^k}{m\neq 0}} f(n-m, n'-Q(n,m))K(m),
\eeq
where $(n, n') \in \Z^k \times \Z^l$, $m\in \Z^k$,  $Q$ is a polynomial mapping from $Z^k \times \Z^k$ to $\Z^l$, and $K$ is a Calder\'{o}n-Zygmund kernel. 
Significant work of \cite{SW1}, \cite{IMSW} has recently proved that $R$, as defined in (\ref{quasi_R}), extends to a bounded operator on $\ell^2(\Z^k \times \Z^l)$ for $Q$ of any degree, and furthermore that $R$ is bounded on  $\ell^p(\Z^k\times \Z^l)$ for all $1<p<\infty$ if  the polynomial mapping $Q$ is of degree  at most 2. 

In particular, this implies that for $\lam = 1+ i\ga$ with $\ga  \neq 0$, the operator $T^\lam$ we study in this paper is bounded on $\ell^p(\Z^{2k+1})$ for all $1< p< \infty$; for $0 < \lam < 1$, the operator $T^\lam$ as defined in (\ref{T_dis}) is a fractional version of (\ref{quasi_R}). Theorem \ref{frxnl_int_thm_dim1} represents the first result for a discrete fractional Radon transform in a quasi-translational invariant setting. 
Most of the recent successes in studying discrete operators have come from intricate decompositions of the operators motivated by the circle method of Hardy and Littlewood, an approach initiated in the setting of discrete operators by Bourgain \cite{Bour88A}.
 In particular, the $\ell^p$ bounds for (\ref{quasi_R}) proved in \cite{SW1} and \cite{IMSW} employed sophisticated techniques rooted in the circle method.

In contrast, the method presented here to prove Theorem \ref{frxnl_int_thm_dim1} is a simple argument with a combinatorial flavor. It is motivated by work of Oberlin \cite{Obe} on $\ell^p \maps \ell^q$ bounds for two translation invariant discrete fractional integral operators acting on functions of $\Z$. Oberlin employed a lemma of Christ on the distribution functions of iterations of operators in order to reduce the problem to counting the number of elements in an arbitrary finite set of integers that are of a certain form (for example, a sum of three squares).
We develop this idea in a new direction to treat the quasi-translation invariant operator $T^\lam$. The limited utility of the Fourier transform for this operator makes this simple method quite attractive. 
The main aspects of the  approach demonstrated here can be generalized to higher dimensions, but the key final step appears to fail for dimensions $k \geq 2$. Ultimately we reduce the problem to bounding the number of integer solutions to a single Diophantine equation; in higher dimensions this equation involves far too many variables, relative to the degrees of freedom allowed in a successful bound for the operator, to admit acceptable bounds.

We proceed in this paper as follows: in Section \ref{sec_reduction} we reduce Theorem \ref{frxnl_int_thm_dim1} to a weak-type bound for an operator acting on characteristic functions of finite sets. In Section \ref{sec_iterate} we perform an iterative procedure, further reducing the problem to bounding a certain finite sum, which we then do in Section \ref{sec_bound} by counting integer solutions to a Diophantine equation. In Section \ref{sec_high_dim} we briefly outline an approach for higher dimensions, and point out the difficulties that arise. Finally, in Section \ref{sec_T_nec} we outline examples that show that conditions (i) and (ii) in Theorem \ref{frxnl_int_thm_dim1} are necessary for $T^\lam$ to be bounded from $\ell^p$ to $\ell^q$.

\section{Reduction of the problem}\label{sec_reduction}
We first reduce Theorem \ref{frxnl_int_thm_dim1} to a restricted weak-type estimate. Note that $\lam=1/3$ is the crossover value at which conditions (i) and (ii) (with equalities) meet at the single point $(1/p, 1/q)=(2/3, 1/3)$. Therefore by interpolation, in order to prove Theorem \ref{frxnl_int_thm_dim1}, it is sufficient to prove that $T^\lam$ is of restricted weak-type $(3/2,3)$ whenever $\lam=1/3 + \ep$, with $\ep >0$. Indeed, complex interpolation with the $(\ell^p,\ell^p)$ bound for $T^\lam$ for $\lam=1+i\ga$ with $\ga \neq 1$, due to \cite{IMSW} (in fact the trivial $(\ell^p,\ell^p)$ bound for $T^\lam$ for any $\Re(\lam)>1$ would suffice), followed by applying the inclusion property of $\ell^p$ spaces (namely $\ell^{q_2} \subset \ell^{q_1}$ if $q_2 < q_1$) and taking adjoints, then proves Theorem \ref{frxnl_int_thm_dim1} for $\lam>1/3$. (Note that a restricted weak-type $(3/2,3)$ bound precisely at $\lam=1/3$ would imply best possible results, namely $(\ell^p,\ell^q)$ bounds for $p,q$ satisfiying equality in condition (i).) For $\lam \leq 1/3$, condition (ii) is stronger than condition (i), and the desired results follow from interpolating the result for $\lam=1/3+\ep$ with the trivial $(\ell^1, \ell^\infty)$ bound for $T^\lam$ for any $\Re(\lam) \geq 0$, followed by applying the inclusion property and taking adjoints.

Thus suppose that $\lam=1/3+\ep$ for a given fixed $\ep>0$. It suffices to consider the action of $T^\lam$ on non-negative functions $f$. We make a dyadic decomposition $T^\lam = \sum_{j=0}^\infty T^\lam_j$, setting
\beq\label{T_dyadicHeis}
 T_j^\lam f(n,t) = 2^{-2\lam j} \sum_{\bstack{m \in \Z^2}{2^j \leq |m| < 2^{j+1}}} f(n-m, t-\om(n,m))
 \eeq
for each $j \geq 0$. Then it is sufficient to show that the dyadic operators $T^\lam_j$ are uniformly weak-type $(3/2,3)$ for all $j \geq 0$.
Fix $j$, and let $T$ denote $T^\lam_j$. By a well-known criterion for restricted weak-type operators (see Chapter 5 of \cite{SteinWeiss}), it is sufficient to prove that for all characteristic functions $\chi_E$ of measurable sets $E$, 
\beq\label{weak-type}
 ||T\chi_E||_{q,\infty}^* \leq A ||E||_{p,1}^* = A|E|^{1/p},
 \eeq
where $(p,q) = (3/2,3)$ and $||\cdot||_{p,1}^*$ and $||\cdot||_{q,\infty}^*$ are weak-type norms.\footnote{
Precisely, the norms are defined as follows (see Chapter 5 of \cite{SteinWeiss}). Let $g^*$ be a non-increasing rearrangement of $g$. Then
$ ||g||^*_{q,\infty} = \sup_{t>0} t^{1/q} g^*(t);$ note that this being finite is equivalent to $\al^q \{ x: |g| >\al\} < A$ for all $\al >0$. Also, 
\[ ||g||_{p,q}^* = \left( \frac{q}{p} \int_0^\infty [t^{1/p}g^*(t)]^{q} \frac{dt}{t} \right)^{1/q}.\]
Conveniently, for a measurable set $E$, $||\chi_E||^*_{p,q} = |E|^{1/p}$ for all $1 \leq q \leq \infty$.}
The key step in proving (\ref{weak-type}) is the following lemma of Christ \cite{Chr98}, which we quote in the form given in \cite{Obe}:
\begin{lemma}[Christ]\label{Christ_lemma}
Suppose that $T$ is a positive operator taking measurable functions to measurable functions and in particular that $T$ has this property for characteristic functions $\chi_E$ of measurable sets $E$. Given $\al>0$ and a measurable set $E$ with $|E|>0$, set 
\begin{eqnarray*}
F & = & \{x: \al < T \chi_E(x) < 2\al\}, \\
\be & = & |E|^{-1}\langle \chi_F, T \chi_E \rangle,
\end{eqnarray*}
when $ |F|>0$.
Then for every $r \geq 0$ there exist constants $\del_r, \ep_r >0$ such that if we define $E_0=E$, $F_0=F$, and set
\begin{eqnarray*}
E_{r+1} &=& \{x \in E_r : T^* \chi_{F_r} (x) \geq \del_r \be \} \\
F_{r+1} &= &\{x \in F_r : T \chi_{E_{r+1}} (x) \geq \ep_r \al \},
\end{eqnarray*}
then all the sets $E_r, F_r$ are nonempty.
\end{lemma}
In Christ's work, this lemma arose in the context of studying the $L^p(\R^n) \maps L^q(\R^n)$ mapping properties of convolution operators defined in terms of a measure supported on the curve $(t, t^2, \ldots, t^n)$. In general, the conclusion that $E_r, F_r$ are nonempty follows from the statements that for sufficiently small $\eta_r >0$, 
 \begin{eqnarray}
 \langle T^* \chi_{F_r}, \chi_{E_{r+1}} \rangle &\geq& \eta_{r} \be |E|, \label{Christ1} \\
\langle \chi_{F_{r+1}}, T \chi_{E_{r+1}} \rangle &\geq& \eta_{r+1} \al |F| \label{Christ2}.
\end{eqnarray}
These statements may be proved by induction. The base cases follow directly from the definitions. Indeed, in the case $r=0$,
\[ \langle T^* \chi_{F}, \chi_{E_{1}}\rangle = \langle \chi_{F}, T \chi_{E} \rangle
	 - \langle T^* \chi_{F}, \chi_{E \setminus E_{1}} \rangle \geq \be |E| - \del_0 \be |E|,\]
where the last term arises because by definition $T^* \chi_{F} < \del_0 \be$ on $E \setminus E_{1}$. Thus choosing $\eta_0 = \del_0 =1/2$, (\ref{Christ1}) holds.
To prove (\ref{Christ2}) in the case $r=0$, note that
\[ \langle  \chi_{F_1}, T\chi_{E_{1}}\rangle = \langle T^* \chi_{F},  \chi_{E_1} \rangle
	 - \langle \chi_{F \setminus F_1}, T\chi_{E} \rangle \geq \eta_0 \be |E| - \ep_0 \al |F|,\]
	 where we have applied (\ref{Christ1}) and the fact that $T^*\chi_F< \ep_0\al$ on $F\setminus F_1$. 
 Recall that $\be = |E|^{-1} \langle \chi_F, T \chi_E \rangle,$ and for a point ${\bf x}$ to be in $F$, one must have $\al < T_{\chi_E} ({\bf x}) < 2\al$, and thus $\be \approx \al |F| |E|^{-1}$. 
	 Thus choosing $\eta_1=\ep_0 = \eta_0/2$, (\ref{Christ2}) follows.

Assuming now that the two relations (\ref{Christ1}) and (\ref{Christ2}) hold for $0, \ldots, r$, then
\[ \langle T^* \chi_{F_r}, \chi_{E_{r+1}}\rangle = \langle \chi_{F_r}, T \chi_{E_r} \rangle
	 - \langle T^* \chi_{F_r}, \chi_{E_r \setminus E_{r+1}} \rangle \geq \eta_{r}\al |F| - \del_r \be |E_r|,\]
where the last term arises because by definition $T^* \chi_{F_r} < \del_r \be$ on $E_r \setminus E_{r+1}$. But since $\be \approx \al |F||E|^{-1}$, in fact
\[ \langle T^* \chi_{F_r}, \chi_{E_{r+1}}\rangle \geq \eta_{r} \be |E|  - \del_r \be |E_r|.\]
Now choosing $\eta_{r+1} = \del_r = \eta_{r}/2$, the relation (\ref{Christ1}) for $r+1$ follows, since $|E_r| \leq |E|$.
 The inductive process for (\ref{Christ2}) proceeds in the same manner, and the lemma follows.

\section{Iteration}\label{sec_iterate}
With the definitions of Lemma \ref{Christ_lemma}, to prove (\ref{weak-type}) it is sufficient to show that $\al^q |F| \leq A |E|^{q/p}$, i.e.
\beq\label{desired_ineq}
\al^3 |F| \leq A |E|^{2}, 
\eeq
for all non-empty measurable sets $E \subset \Z^{2} \times \Z$ (using the counting measure).
Note that for $T=T^\lam_j$, the dual operator $T^*$ is given by
\[ T^*g(n,t) = 2^{-2\lam j} \sum_{\bstack{m \in \Z^2} {2^j\leq |m| < 2^{j+1}}} g(n+m, t+\om(n,m)).\]
Therefore, given $r \geq 1$, choose some ${\bf{x}}=(x,x') \in E_r \subseteq \Z^2 \times \Z$, which exists by Christ's lemma. By the definition of $E_r$, this implies $T^* \chi_{F_{r-1}} (x, x') \geq \del_{r-1} \be$, i.e.
\[ \sum_{|m| \approx 2^j} \chi_{F_{r-1}} (x+m, x'+ \om(x,m)) \geq 2^{2\lam j} \del_{r-1} \be,\]
or, using the notation $\Om^*({\bf x};m) = (x+m, x'+ \om(x,m))$,
\beq\label{E_r}
 \sum_{|m| \approx 2^j} \chi_{F_{r-1}} (\Om^*({\bf x};m)) \geq 2^{2\lam j} \del_{r-1} \be.
 \eeq
Similarly, if ${\bf y} = (y,y') \in F_r$ then $T \chi_{E_r} (y,y') \geq \ep_{r-1} \al$, i.e.
\beq\label{F_r}
 \sum_{|m| \approx 2^j} \chi_{E_r} (\Om({\bf y};m)) \geq 2^{2\lam j} \ep_{r-1} \al,
 \eeq
where $\Om({\bf{y}};m) = (y-m, y'-\om(y,m))$. Note that both $\Om, \Om^* : \Z^3 \times \Z^2 \maps \Z^3$.

Starting with a point in a certain set $E_r$ or $F_r$, the inequalities (\ref{E_r}), (\ref{F_r}) allow one to step backwards sequentially to the given arbitrary set $E$, each time counting how many points remain. The set in which we start this process is determined by the exponents in the final inequality (\ref{desired_ineq}) that is our goal. In our case, we fix ${\bf y} \in F_1$ and consider the expression 
\beq\label{S_expression}
 S = \sum_{\bstack{m_1, m_2, m_3 \in \Z^2}{|m_l| \approx 2^j}} \chi_E(\Om (\Om^*(\Om({\bf y};m_1);m_2);m_3)) .
 \eeq
We now proceed to give a lower bound for the sum $S$ in terms of $|F|$ and $|E|$.

By (\ref{F_r}), there exist at least $2^{2\lam j}\ep_0 \al$ values of $m_1$ with $\Om({\bf y}; m_1) \in E_1$.  By (\ref{E_r}), there then exist at least $2^{2\lam j}\del_0 \be$ values of $m_2$ with $\Om^*(\Om({\bf y};m_1);m_2) \in F_0.$ Now recall that $F_0 = F$ and $F=\{ x: \al < T \chi_E < 2\al\}$, so if ${\bf y} \in F$
then $T \chi_E( {\bf y}) = 2^{-2\lam j} \sum_{|m| \approx 2^j} \chi_E(\Om({\bf y};m)) > \al$; thus there exist at least $2^{2\lam j} \al $ values of $m_3$ such that $\Om(\Om^*(\Om({\bf y};m_1);m_2);m_3)$ lies in $E_0 = E$. 
Assembled in sequence, these steps prove the lower bound
\[ S \geq 2^{6 \lam j} \ep_0 \del_0 \al^2 \be.\]

Recall that $\lam = 1/3+\ep$; let $\del = 3\ep$, so that $3 \lam = 1 + \del$; also recall that $\be \approx \al |F|\; |E|^{-1}$.
Therefore, for some constant $c$ (dependent on $\ep_0, \del_0$, which are fixed), 
\beq\label{S_lowerbd}
 S \geq c2^{(2 + 2\del) j} \al^3 |F| |E|^{-1}.
 \eeq
Thus to show the desired inequality (\ref{desired_ineq}), it is now sufficient to show that 
\beq\label{S_ineq}
S \leq c 2^{(2+2\del)j} |E|.
\eeq

\section{Bounding $S$}\label{sec_bound}
Note that by definition
\begin{eqnarray*}
 \Om(\Om^*(\Om({\bf y};m_1);m_2);m_3)
 & = &  \Om(\Om^*(y-m_1, y'-\om(y,m_1);m_2);m_3) \\
 & = & \Om(y-m_1+m_2, y' - \om(y,m_1) + \om(y-m_1,m_2);m_3) \\
 & = &  (y-m_1+m_2 - m_3, y' - \om(y,m_1) + \om(y-m_1,m_2) \\
 	&& - \; \om(y-m_1 + m_2, m_3)).
 \end{eqnarray*}
In order to give an upper bound for $S$, we need to count how many points of this shape are in $E$; replacing the set $E$ by the set ${\bf y}-E$, which we rename $E$, this is equivalent to counting the number of points in $E$ of the form 
\[(m_1-m_2 + m_3,  \om(y,m_1) - \om(y-m_1,m_2) + \om(y-m_1 + m_2, m_3)).\]

Therefore we now define 
\[ \tilde{S} = \sum_{(a,b,c) \in E} \; \; \sum_{m_1, m_2, m_3} 1,\]
where the inner sum is restricted to those $m_1, m_2, m_3 \in \Z^2$ with $|m_l| \approx 2^{j}$ that
satisfy the conditions 
\begin{eqnarray}
a& =& m_{1,1} - m_{2,1} + m_{3,1}\label{cond_1}\\
b& = & m_{1,2} - m_{2,2} + m_{3,2} \label{cond_2} \\
c&=& \om(y,m_1) - \om(y-m_1,m_2) + \om(y-m_1 + m_2, m_3) . \label{cond_3}
\end{eqnarray}
To prove (\ref{S_ineq}) it is now sufficient to show 
\beq\label{S_tilde_ineq}
\tilde{S} \leq c2^{(2+2\del)j} |E|.
\eeq
Given $(a,b,c) \in E$, set $A = (a,b)$. Define $n_1, n_2, n_3 \in \Z^2$ by 
\beq\label{m_n_dfn}
 m_1 = \frac{n_1 + A}{3}, \quad m_2 = \frac{n_2 - A}{3}, \quad m_3 = \frac{n_3 + A}{3}.
 \eeq
Then conditions (\ref{cond_1}), (\ref{cond_2})  become (in vector form)
\[ \frac{n_1 + A}{3} - \frac{n_2 - A}{3} + \frac{n_3 + A}{3} = A,\]
or 
\beq\label{cond_1n}
 n_1 - n_2 + n_3 =0.
\eeq
Condition (\ref{cond_3}) evolves under the change of variables from $m_1, m_2, m_3$ to $n_1, n_2, n_3$ into 
\begin{multline*}
 3\om(y,n_1) - 3\om(y,n_2) + 3\om(y,n_3) + \om(n_1+A, n_2 - A) \\
   - \om(n_1 + A, n_3 + A)
	+ \om(n_2 -A, n_3 + A)  = c_1,
	\end{multline*}
where several constant terms such as $\om(y,A)$ and $\om(A,A)$ have been moved to the right hand side. After further simplification, this becomes
\begin{multline*}
  3\om(y,n_1) - 3\om(y,n_2) + 3\om(y,n_3)
	 + \om(n_1,n_2) + \om(n_2,n_3) \\
	   - \om(n_1, n_3) -2\om(n_1,A) - 2\om(A, n_3) = c_2.
\end{multline*}
We now use (\ref{cond_1n}) to eliminate $n_2 = n_1+n_3$, which results in
\[ \om(n_1,n_3) - 2\om(n_1, A) - 2\om(A, n_3) = c_2.\] 
Our goal is to write this in a form that allows us to restrict as many variables as possible to a negligible number of choices; note that it is equivalent to the condition 
\beq\label{om_next}
 \om(n_1,n_3-2A) +\om(-2A, n_3) = c_2.
 \eeq
But since the form $\om$ is symplectic, $\om(-2A,-2A)=0$, so (\ref{om_next}) is equivalent to 
\beq\label{symp_form}
  \om(n_1,n_3-2A) +\om(-2A, n_3-2A) =  \om(n_1-2A,n_3-2A)  = c_2.
  \eeq
We now define $t_1, t_2 \in \Z^2$ by $t_1 = n_1-2A$, $t_2 = n_3 - 2A$, so that (\ref{symp_form}) is equivalent to
\[ \om(t_1,t_2)= 2t_{1,2}t_{2,1} - 2t_{1,1}t_{2,2}= c_2.\]

This is an equation in four variables; choosing $t_{1,1}$ and $t_{2,2}$ freely (for which there are $2^{2j}$ choices), the equation becomes of the form $t_{1,2}t_{2,1} = c_3$. Since $t_{1,2}, t_{2,1}\approx 2^j$, for this to have any solutions we must have $c_3 \approx 2^{2j}$. Furthermore, the number of choices for the remaining two variables $t_{1,2}, t_{2,1}$ is limited by the number of divisors of $c_3$, which is $O(c_3^\eta) = O(2^{2\eta j})$ for every $\eta >0$ (see for example Theorem 315 of  \cite{HardyWright}).

There are therefore  $O(2^{2j + 2\eta j})$ total choices for the (vector) variables $t_1, t_2$, and since these specify $n_1, n_3$ uniquely, there are $O(2^{2j + 2\eta j})$ choices for $n_1, n_3$. Note that condition (\ref{cond_1n}) specifies $n_2$ uniquely once $n_1, n_3$ are chosen, hence there are  $O(2^{2j + 2\eta j})$ choices for $n_1, n_2, n_3$. Finally, (\ref{m_n_dfn}) defines $m_1,m_2,m_3$ uniquely in terms of $n_1, n_2, n_3$, thus there are  $O(2^{2j + 2\eta j})$ choices for $m_1, m_2, m_3$, given any element $(a,b,c) \in E$. This proves the bound (\ref{S_tilde_ineq}) for $\tilde{S}$, and hence the final weak-type inequality (\ref{desired_ineq}), from which Theorem \ref{frxnl_int_thm_dim1} follows immediately.

\section{Higher dimensions}\label{sec_high_dim}
In higher dimensions $k \geq 2$,
we expect that for $0<\lam<1$, the discrete operator $T^\lam$ defined in (\ref{T_dis}) extends to a bounded operator from $\ell^p(\Z^{2k+1})$ to $\ell^q(\Z^{2k+1})$ if and only if $p,q$ satisfy:
\begin{enumerate}
\item  $1/q \leq 1/p - \frac{2k(1-\lam)}{2k+2}$,
\item $1/q< \lam, 1/p> 1-\lam$.
\end{enumerate}
Predictably, the method described for the case $k=1$ becomes much more complicated in higher dimensions, because the chain of $E_r, F_r$ deductions grows with the dimension. Nevertheless, it is possible to extend the main steps of the method to arbitrary dimensions, but it appears that the final step of providing a good upper bound for integer solutions to the higher-dimensional analogues of equations (\ref{cond_1}) -- (\ref{cond_3}) is not feasible. We outline a general argument and briefly summarize the difficulties below. 

In general, in order to show that $T^\lam$ maps $\ell^p(\Z^{2k+1}) \maps \ell^q(\Z^{2k+1})$ for $\lam, p,q$ satisfying the requirements (i) and (ii) but with a strict inequality in (i), it is sufficient to prove a restricted weak-type estimate just above the crossover value $\lam = 1/(k+2)$ at which conditions (i) and (ii) (with equalities)  meet at the single point $(1/p,1/q) = (\frac{k+1}{k+2},\frac{1}{k+2})$. 
Therefore, for $\lam=1/(k+2) + \ep$, with any fixed $\ep>0$, we make a dyadic decomposition of the operator $T^\lam = \sum_{j=0}^\infty T_j^\lam$ as in (\ref{T_dyadicHeis}):
\[ T_j^\lam f(n,t) = 2^{-2k\lam j} \sum_{\bstack{m \in \Z^{2k}}{2^j \leq |m| < 2^{j+1}}} f(n-m, t-\om(n,m)).\]
It then suffices to prove that each $T_j^\lam$ is uniformly of weak-type $(\frac{k+2}{k+1}, k+2)$. From now on, we fix $j$ and refer to $T_j^\lam$ as $T$. As in the case $k=1$, we may restrict our attention to the action of $T$ on characteristic functions of measurable sets $E \subset \Z^{2k+1}$, and it suffices to prove the analogue of (\ref{weak-type}) for $(p,q)=(\frac{k+2}{k+1}, k+2)$.
Again defining $F = \{ x: \al < T_{\chi_E}(x) < 2\al\}$, the problem reduces to showing that for any $\al>0$, and any measurable set $E$,  $\al^q |F| \leq A |E|^{q/p}$, which in this case is the inequality
\beq\label{FE_ineq}
 \al^{k+2} |F| \leq A |E|^{k+1}.
 \eeq

For notational convenience, we define two auxiliary functions $\Om, \Om^* : \Z^{2k+1} \times \Z^{2k} \maps \Z^{2k+1}$, acting on $\x  = (x, x') \in \Z^{2k} \times \Z$, and $m \in \Z^{2k}$ by 
\begin{eqnarray*}
\Om(\x;m) &=& (x-m, x'-\om(x,m)), \\
\Om^*(\x;m) &=& (x+m, x'+\om(x,m)).
\end{eqnarray*}
Then Christ's lemma shows that for $\x \in E_r$,
\beq\label{E_k'}
 \sum_{|m| \approx 2^j} \chi_{F_{r-1}} (\Om^*({\bf x};m)) \geq 2^{2k\lam j} \del_{r-1} \be,
 \eeq
and similarly for ${\bf y} \in F_r$,
\beq\label{F_k'}
 \sum_{|m| \approx 2^j} \chi_{E_r} (\Om({\bf y};m)) \geq 2^{2k\lam j} \ep_{r-1} \al.
 \eeq

We would like to proceed as in the case $k=1$ by starting with a point $\x \in E_r$ or $\y \in F_r$, for an appropriate $r \geq 1$, and then studying a sum analogous to (\ref{S_expression}), but with an appropriate number of iterations of $\Om, \Om^*$ depending on the dimension $k$. (Since $k+2$ factors of $\al$ appear in the desired inequality (\ref{FE_ineq}), $k+2$ factors of $\Om$ or $\Om^*$ will appear in the chain defining $S$.) 

The difficulties inherent in the method arise already when $k=2$, thus we temporarily restrict our attention to this case. The desired inequality (\ref{FE_ineq}) now takes the form
\beq\label{FE_ineq2}
 \al^{4} |F| \leq A |E|^{3}
 \eeq
 and the relevant sum $S$ is defined by
 \beq\label{S_expression2}
 S=\sum_{\bstack{m_1, m_2, m_3, m_4 \in \Z^4}{|m_l| \approx 2^j}} \chi_E(\Om (\Om^*(\Om(\Om^*(\x;m_1);m_2);m_3);m_4)) .
 \eeq
 Working sequentially downward from a point $\x \in E_2$ via the inequalities (\ref{E_k'}) and (\ref{F_k'}), we see there are at least $2^{2k\lam j}\del_1 \be$ values $m_1$ with $\Om^*(\x;m_1) \in F_1$, at least $2^{2k\lam j} \ep_0 \al$ values $m_2$ with $\Om(\Om^*(\x;m_1);m_2) \in E_1$, at least $2^{2k\lam j} \del_0 \be$ values $m_3$ with $\Om^*(\Om(\Om^*(\x;m_1);m_2);m_3) \in F_0=F$, and finally at least $2^{2k\lam j} \al$ values $m_4$ with $\Om(\Om^*(\Om(\Om^*(\x;m_1);m_2);m_3);m_4) \in E_0=E$. Thus in total,
 \beq\label{S_lowerbd2}
  S \geq c 2^{4 \cdot 2k \lam j} \al^2 \be^2 \geq c2^{16\lam j}\al^4 |F|^2|E|^{-2},
  \eeq
 where again we use  $\be \approx \al |F||E|^{-1}$. 
 Thus in order to prove (\ref{FE_ineq2}), we would need to prove $ S \leq c2^{16\lam j} |E| \; |F|.$ Recalling that $\lam=1/4 +\ep$, we set $\del=4\ep$ so that our goal is to prove that
 \beq\label{SE_EF}
 S \leq c2^{(4+4\del) j} |E| \; |F|.
 \eeq

This is more complicated than the analogous inequality (\ref{S_ineq}) for $k=1$, as it requires comparing $S$ both to $|E|$ and to $|F|$. In fact it would suffice to consider 
 the modified sum $\tilde{S}$ defined by (\ref{S_expression2}), but with the additional requirement that $m_1, m_2,m_3$ be such that $\Om^*(\Om (\Om^*({\bf{x}};m_1);m_2);m_3) \in F.$ This still satisfies
 \[ \tilde{S} \geq c2^{16 \lam j}\al^4 |F|^2 |E|^{-2},\]
so to prove (\ref{FE_ineq2}), it would be sufficient to show that 
\beq\label{SE_EF_tilde}
 \tilde{S} \leq c2^{(4 + 4\del) j} |E| \; |F|.
  \eeq
By definition, $\tilde{S} \leq N_1(F) N_2(E),$ where
\begin{eqnarray*}
N_1(F) & = &   \#\{m_1, m_2, m_3 \in \Z^4 : (x+m_1-m_2+m_3, \\
	&& \quad x' + \om(x,m_1) - \om(x+m_1,m_2) +\om(x+m_1-m_2,m_3)\in F \} \\
N_2(E) & = &  \sup_{(v,v') \in F} N_2(v,v'),
\end{eqnarray*}
where
\[ N_2(v,v')= \#\{m_4 \in \Z^4 : (x+v -m_4, x' + v'  - \om(x+v,m_4) \in E \} .\]
After shifting the set $F$ by the fixed shift ${\bf x}-F$ (which we rename $F$),  
\[ N_1(F) = \sum_{(A,b) \in F} \sum_{m_1,m_2,m_3} 1\]
where  $(A,b) \in F \subset \Z^{4}\times \Z$ and the inner sum is restricted to $m_1,m_2, m_3 \in \Z^4$ with $|m_l| \approx 2^j$ and satisfying the four equations given (in vector form) by
\beq\label{m_relation0}
m_1 -m_2 +m_3 = A 
\eeq
and the fifth equation
\beq\label{om_relation0}
 \om(x,m_1) 
		- \om(x+m_1,m_2) +\om(x+m_1-m_2,m_3) =  b.
		\eeq
Define $n_1, n_2, n_3\in \Z^4$ by 
\[ m_1 = \frac{n_1 + A}{3}, \quad m_2 = \frac{n_2 - A}{3} \quad m_3= \frac{n_3 +A}{3},\] so that (\ref{m_relation0}) becomes 
the condition $n_1-n_2+n_3=0$. Then by the same procedure that led to (\ref{symp_form}), the relation (\ref{om_relation0}) becomes
\[\om(n_1-2A,n_3-2A)  = c\]
for some constant $c$. Defining $t_1,t_2 \in \Z^4$ by $t_1 = n_1-2A$, $n_3 -2A$, it suffices to count solutions to
\[ \om(t_1,t_2) = c,\]
with the coordinates of $t_1,t_2$ all of size $\approx 2^j$. This is an equation in 8 variables;
choosing $6$ coordinates freely (for $O(2^{6j})$ total choices), the last two coordinates are divisors of a constant, and hence there are only $O(2^{\eta j})$ choices for them. Thus $N_1(F) \leq c2^{(6 +\eta)j}|F|$. 

It is trivial to bound $N_2(v,v')$ for each $\mathbf{v} \in F$, since after shifting $E$ by the fixed shift $\mathbf{x} + \mathbf{v}-E$, 
\[ N_2(v,v') = \sum_{(A,b) \in E} \sum_{m_4} 1,\]
where the sum is restricted to those $m_4$ such that $m_4=A$ and $\om(x+v,m_4) = b$. There is only one choice of $m_4$ for each $A$, hence $N_2(v,v') \leq |E|$ and $N_2(E) \leq |E|$. Along with the bound for $N_1(F)$, this ultimately leads to the bound
\[ \tilde{S} \leq c 2^{(6+2\eta)j}|F| \; |E| \]
for any $\eta>0$. Unfortunately, this is too large for the desired bound (\ref{SE_EF_tilde}).

In general, the $k$-dimensional analogue to $S$ can be handled similarly. It includes a chain of $k+2$ applications of $\Om$ or $\Om^*$ and it can be broken up into as many separate counting functions like $N_1(F)$ and $N_2(E)$ as necessary to match the powers of the cardinalities of $|E|$ and $|F|$ in the analogue to (\ref{SE_EF}). But the resulting upper bounds are not sufficiently sharp: the $k$-dimensional analogues of equations (\ref{cond_1}) -- (\ref{cond_3}) take the form of $2k+1$ equations in $2k(k+2)$ variables, which after elimination of variables via the analogue of (\ref{cond_1n}) leads to 1 equation in $2k(k+1)$ variables. But the analogue to (\ref{S_ineq}) only allows $2k$ degrees of freedom, and thus it appears that it is not possible to obtain a sufficiently sharp bound to prove the desired relation $\al^{k+2}|F| \leq A |E|^{k+1}$ in higher dimensions.

\section{Necessary conditions}\label{sec_T_nec}
In this section we indicate why the conditions (i), (ii) in Theorem \ref{frxnl_int_thm_dim1} are in fact necessary for $T^\lam$ to map $\ell^p$ to $\ell^q$. As it will not cause any difficulty to consider higher dimensions in these examples, we show that the conditions 
\begin{enumerate}
\item  $1/q \leq 1/p - \frac{2k(1-\lam)}{2k+2}$
\item $1/q< \lam, 1/p> 1-\lam$
\end{enumerate}
are necessary for the operator $T^\lam$ defined in (\ref{T_dis}) to map $\ell^p(\Z^{2k+1})$ to $\ell^q(\Z^{2k+1})$, in any dimension $k \geq 1$.

For (ii), let $f(n,t) = 1$ if $(n,t)=0$ and let $f$ vanish otherwise, so that clearly $f \in \ell^p$ for every $1 \leq p \leq \infty$. One then sees that $T^\lam f(n,t) = |n|^{-2k\lam}$ if $t=0$ and vanishes otherwise, since the condition $n-m=0$ forces $\om(n,m) = \om(n,n)$ to vanish, and as a result $t$ must be zero if $t-\om(n,m) =0$. Thus
\[ ||T^\lam f||^q_{\ell^q} = \sum_{n \in \Z^{2k}} |T^\lam f (n,0)|^q = \sum_{n \in \Z^{2k}} |n|^{-2kq\lam},\]
which is finite if $1/q<\lam$. The condition $1/p >1-\lam$ follows by taking adjoints.

For the diagonal condition (i) we must work a bit harder. 
Define $f(n,t) = |t|^{-\al} \chi(n/|t|^{1/2})$ for some $\al >0$ to be chosen later, where $\chi$ denotes the characteristic function of the ``box'' $\{x \in \R^{2k}: 1/2 < |x_j| <2, 1 \leq j \leq 2k\}$. Then 
\[ ||f||_{\ell^p}^p = \sum_{n,t} |f(n,t)|^p = \sum_t \frac{1}{|t|^{\al p}} \sum_{\bstack{n_1, \ldots, n_{2k}}{\frac{1}{2}|t|^{1/2} \leq |n_j| \leq 2|t|^{1/2}}} \chi(\frac{n}{|t|^{1/2}})^p \approx c \sum_t \frac{|t|^{2k/2}}{|t|^{\al p}}. \]
This is finite if $\al > \frac{k+1}{p}$.
Now consider 
\begin{multline}\label{T23}
 T^\lam f(n,t) = \sum_{\bstack{m \in \Z^{2k}}{m \neq 0}} \frac{f(n-m, t-\om(n,m))}{|m|^{2k\lam}} \\
  \geq \sum_{\bstack{m}{|\om(n,m)| \leq \del |t|}} \frac{1}{|t-\om(n,m)|^\al} \chi \left( \frac{n-m}{|t -\om(n,m)|^{1/2}} \right) \frac{1}{|m|^{2k\lam}},
 \end{multline}
 for some small $0< \del < 1/2$. 
We can impose the restriction on $m$ that $|\om(n,m)| \leq \del |t|$ and obtain a lower bound for the operator value because all the summands are non-negative; furthermore, under this assumption, $t-\om(n,m) \approx t$. In fact, we can further restrict $m$ so that $|m| < \del |t|^{1/2}$, in which case 
\beq\label{frac_eqn}
 T^\lam f(n,t) \geq \frac{1}{|t|^{\al}} \chi(\frac{n}{|t|^{1/2}}) \sum_{|m| \leq \del |t|^{1/2}} \frac{1}{|m|^{2k\lam}} =  \frac{1}{|t|^{\al}} \chi (\frac{n}{|t|^{1/2}} ) \sum_{l=1}^{\del^2 |t|} \frac{r_{2k}(l)}{l^{k\lam}}  ,
  \eeq
where $r_{2k}(l)$ denotes the number of representations of $l$ as a sum of $2k$ squares.
On average,
\[ \sum_{l=1}^L r_{2k}(l) = \frac{\pi^k}{\Ga(k+1)}L^{k} + o(L^k);\]
 this may be shown by comparing the number of integer lattice points in the Euclidean ball of radius $l^{1/2}$ to the volume of the ball (see \cite{Wal} for example).
Thus by partial summation, the 
 last sum in (\ref{frac_eqn}) is $\approx |t|^{k(1-\lam)}$ (plus a smaller error term), and hence
\[ T^\lam f(n,t) \geq C \frac{1}{|t|^\al} \chi(\frac{n}{|t|^{1/2}})|t|^{k(1-\lam)}.\]
Therefore, 
\[
 ||T^\lam f||^q_{\ell^q} \geq C\sum_t |t|^{kq(1-\lam)-\al q} \sum_n \chi (\frac{n}{|t|^{1/2}})^q 
 \approx C\sum_t |t|^{kq(1-\lam) -\al q + k} ,\]
 since there are about $|t|^{2k/2} = |t|^k$ points with $|n| \approx |t|^{1/2}$.
This last sum is finite if and only if $\al - \frac{k+1}{q} > k(1-\lam).$ Recall that for $f$ to be in $\ell^p$ we required $\al>\frac{k+1}{p}$. Thus set $\al = \frac{k+1}{p} + \ep$ for any $\ep >0$. Then to have $T^\lam f \in \ell^q$, we must have $ \frac{k+1}{p} - \frac{k+1}{q} > k(1-\lam) -\ep$ for any $\ep>0$, and hence $1/p - 1/q \geq \frac{k(1-\lam)}{k+1}.$ This proves the necessity of condition (i).

\section*{Acknowledgements}
 The author would like to thank Elias M. Stein for suggesting this area of inquiry and for his generous advice and encouragement. This work was supported by the Princeton University Centennial Fellowship; the Simonyi Fund at the Institute for Advanced Study; and the National Science Foundation [DMS-0902658, DMS-0635607].

\bibliographystyle{amsplain}
\bibliography{AnalysisBibliography}

\end{document}